# The Fallacy in the Paradox of Achilles and the Tortoise

James Q. Feng

james.q.feng@gmail.com

**Abstract**

Zeno's ancient paradox depicts a race between swift Achilles and a slow tortoise with a head start. Zeno argued that Achilles could never overtake the tortoise, as at each step Achilles arrived at the tortoise's former position, the tortoise had already moved ahead. Though Zeno's premise is valid, his conclusion that Achilles can "never" pass the tortoise relies on equating infinite steps with an infinite amount of time. By modeling the sequence of events in terms of a converging geometric series, this paper shows that such an infinite number of events sum up to a finite distance traversed in finite time. The paradox stems from confusion between an infinite number of events, which can happen in a finite time interval, and an infinite amount of time. The fallacy is clarified by recognizing that the infinite number of events can be crammed into a finite time interval. At a given speed difference after a finite amount of time, Achilles will have completed the infinite series of gaps at the "catch-up time" and passed the tortoise. Hence this paradox of Achilles and the tortoise can be resolved by simply adding "before the catch-up time" to the concluding statement of "Achilles would never overtake the tortoise".

---

One of the most famous paradoxes is that of Zeno's Achilles and the tortoise (for a historical account, cf. Lynds 2003; Theunissen & Oud 2021). An ancient Greek philosopher, Zeno of Elea, imagined a race between the mythical warrior Achilles and a tortoise. To make it seem like a fair race, Achilles allows the tortoise to start at a distance ahead of him. When Achilles reaches the point where the tortoise started, the tortoise has moved ahead of him, and each time Achilles reaches the tortoise's former position, the tortoise is further down the track. Thus, the tortoise is always ahead of Achilles, which leads to the conclusion that Achilles would never be able to run past the naturally slow tortoise—obviously against the common-sense-based expectation.

Zeno's premise is that "For Achilles to catch the tortoise, he must run to the tortoise's initial place when the tortoise would have advanced to a new place, to where Achilles must run to reach when the tortoise would have advanced further ahead to a newer place, and so on." There is nothing wrong with Zeno's premise, which can be mathematically expressed as follows.

If Achilles starts at $x = 0$, running at a constant speed of $s_A$, while the tortoise starts at $x = x_0$, moving at a constant speed of $s_T$, he should reach $x = x_0$ at time $t_0 = x_0 / s_A$ when the tortoise would have moved to a new place at $x = x_1 = x_0 + s_T\, t_0 = x_0\, (1 + s_T / s_A)$ further ahead. Then, Achilles should reach $x = x_1$ at time $t_1 = x_1 / s_A = (x_0 / s_A)\, (1 + s_T / s_A)$, when the tortoise would have moved to a newer place at $x = x_2 = x_0 + s_T t_1 = x_0 [1 + (s_T / s_A) (1 + s_T / s_A)]$ further ahead, and so on. Thus, Achilles should reach $x = x_n$ at time

$$t_n = (x_0 / s_A) * [1 + (s_T / s_A) + (s_T / s_A)^2 + \ldots + (s_T / s_A)^n] = (x_0 / s_A) [1 - (s_T / s_A)^{n+1}] / (1 - s_T / s_A) \quad (1)$$

as a geometric series, with

$$x_n = x_0 [1 - (s_T / s_A)^{n+1}] / (1 - s_T / s_A), \quad (2)$$

where the speed of Achilles and the speed of the tortoise,

$$s_A = x_n / t_n \text{ and } s_T = (x_{n+1} - x_0) / t_n$$

are constants.

Given that the speed ratio $s_T / s_A < 1$, as Achilles is faster than the tortoise, we have

$$t_\infty = (x_0 / s_A) / (1 - s_T / s_A) \text{ and } x_\infty = x_0 / (1 - s_T / s_A) \text{ as } n \to \infty, \quad (3)$$

which suggests that Achilles should catch up with the tortoise at $x = x_\infty$ (the catch-up distance), when $t = t_\infty$ (the catch-up time); thereafter, Achilles would be running ahead of the tortoise. For example, if $x_0 = 1$ and $s_T / s_A = 1/2$ (or 1/3, 1/5, 1/10), assuming $s_T = 1$, we should have the catch-up time $t_\infty = 1$ (or 1/2, 1/4, 1/9) and catch-up distance $x_\infty = 2$ (or 3/2, 5/4, 10/9).

Although Zeno presented his premise flawlessly – i.e., the tortoise would be ahead of Achilles at any time $t = t_n < t_\infty$, no matter how large the value of $n$ becomes – his conclusion that Achilles would *never* be able to run past the tortoise is mistakenly based on an implication that $n \to \infty$ equates to $t_n \to \infty$, means that something could not happen at any time. Yet, "unbounded $n$ corresponds to unbounded $t_n$" was not included in Zeno's premise. In fact, a large $n$ in $t_n$ cannot mean that the value of $t_n$ must be large. With $s_T / s_A < 1$, $t_n < t_\infty = (x_0 / s_A) / (1 - s_T / s_A)$ is always finite, beyond which (as the clock ticks and time moves forward when $t > t_\infty$) Achilles will be ahead of the tortoise.

The paradox of Zeno's Achilles and the tortoise has been regarded as a supertask, consisting of an infinite sequence of subtasks to be completed in a finite amount of time. Whether a supertask can possibly be completed has been a subject of academic debate (e.g., Black 1951; Wisdom 1952; Chihara 1965; McLaughlin and Miller 1992; Alper and Bridger 1997; Salmon 1998; Peijnenburg and Atkinson 2008; Ardourel 2015). The difficulty appears to come from the intuition that an infinite number of subtasks could take forever to complete. But in the case of Zeno's Achilles and the tortoise, the completion of those infinite steps seems possible when one recognizes the fact that an infinite number of subtasks may not necessarily need an infinite amount of time to complete, because the finite catch-up time $t_\infty$ can be finite.

In philosophy, the concepts of space and time have been related to the relative locations of objects and sequences of events. It has become common sense to measure the distance between objects with yardsticks of constant length and to create units of time using a time interval between two specific and regular events, such as repeated positions of the Sun in the sky throughout the day. In some ways, Zeno tricked us by implying the time measurement with events of a shrinking time interval, presenting an illusion of an infinite sequence of events occurring through an infinitely long time. As time is traditionally presented as a dense real number axis, a finite time interval contains an infinite number of real numbers to offer opportunities for creating illusions with an infinite sequence of events crammed into a finite small time interval. Another similar paradox of Zeno describes the runner Achilles starting at the starting line of a track and running past half of the distance to the finish line. He then runs past half of the remaining distance and continues in this way over and over without end, getting ever closer to the finish line but never able to reach it. It has the same type of logical fallacy as that in Achilles and the tortoise; so does Whitrow's version of the paradox of a bouncing ball (Whitrow 1980).

It has been suggested that many problems in philosophy arise from misunderstandings about what everyday words actually mean (Wittgenstein 1953). If correctly stated using appropriate language, the paradox of Achilles and the tortoise can be resolved by simply adding "before the catch-up time $t_\infty$ defined in (3)" to the concluding statement of "Achilles would never be able to run past the tortoise". As might be noted, the resolution presented here involves only consistency in logic and mathematics assuming time to be represented as a real-number with infinite divisibility. In the physical world, modern understanding would indicate (Hilbert 1983): "The sort of divisibility needed to realize the infinitely small is nowhere to be found in reality. The infinite divisibility of a continuum is an operation which exists only in thought." Thus discrete solutions have also been investigated (e.g., Ardourel 2015; Theunissen & Oud 2021). But for resolving the paradox *per se* it should be sufficient to focus just on the logical reasoning according to Zeno's conception, as shown herewith.